# Couple microscale periodic patches to simulate macroscale emergent dynamics

Hammad Alotaibi[*]    Barry Cox[†]    A. J. Roberts[‡]

March 3, 2017

## Abstract

This article proposes a new way to construct computationally efficient 'wrappers' around fine scale, microscopic, detailed descriptions of dynamical systems, such as molecular dynamics, to make predictions at the macroscale 'continuum' level. It is often significantly easier to code a microscale simulator with periodicity: so the challenge addressed here is to develop a scheme that uses only a given periodic microscale simulator; specifically, one for atomistic dynamics. Numerical simulations show that applying a suitable proportional controller within 'action regions' of a patch of atomistic simulation effectively predicts the macroscale transport of heat. Theoretical analysis establishes that such an approach will generally be effective and efficient, and also determines good values for the strength of the proportional controller. This work has the potential to empower systematic analysis and understanding at a macroscopic system level when only a given microscale simulator is available.

# Contents

[*]School of Mathematical Sciences, University of Adelaide, South Australia. mailto:hammad.alotaibi@adelaide.edu.au ORCID:0000-0002-4798-7119

[†]School of Mathematical Sciences, University of Adelaide, South Australia. mailto:barry.cox@adelaide.edu.au ORCID:0000-0002-0662-7037

[‡]School of Mathematical Sciences, University of Adelaide, South Australia. mailto:anthony.roberts@adelaide.edu.au ORCID:0000-0001-8930-1552

# 1 Introduction

Computational molecular simulations has become a valuable tool to the point where computation now stand alongside theoretical and experimental methods in addressing problems in materials science. However, the high computational cost often constrains simulations to limited space-time domains (Dove 2008, e.g.). The Equation-Free Multiscale Scheme aims to use such microscale molecular simulations to efficiently compute and predict large macroscale space-time dynamics (Kevrekidis and Samaey 2009; Liu et al. 2015, e.g.). This article focuses on establishing the basis for a novel design of the equation-free scheme in predicting emergent macroscale properties over large space scales by computing atomistic dynamics only on relatively small widely distributed patches (Samaey, Kevrekidis, and Roose 2005, e.g.). In the scenario where a user has coded a microscale simulator with microscale periodicity, our innovation is to show how to couple such small periodic patches so that the overall scheme predicts the correct macroscale spatial dynamics.

Others have previously used an Equation-Free Approach to aid in molecular simulations (Frederix et al. 2007; Roose et al. 2009, e.g.). However, they concentrated on issues associated with long-time integration, whereas here we focus upon designing effective algorithms for large space domains. Future development of a full multiscale equation-free scheme would combine both aspects.

Alternative multiscale methods that have been proposed are based upon analogous simulations at the microscale level. Foe example, in the flow through a porous medium, Hassard et al. (2016) used smoothed particle hydrodynamics on the microscale to estimate macroscale volume averaged fluxes, with a view to forming a two-scale model that appears like a finite volume scheme on the macroscale. For general gradient driven transport processes, Carr, Perré, and Turner (2016) correspondingly proposed an Extended Distributed Microstructure Model where the macroscale flux is determined as the average of microscale fluxes within micro-cells. Both approaches suggest that microscale simulations can be coupled usefully across macroscales. We similarly propose coupling expressed in terms of macroscale quantities (here the temperature). In principle one or more patches may use the same coupling to couple with surrounding continuum simulations. Thus our approach may readily form part of a hybrid molecular-continuum method (Kalweit and Drikakis 2011, e.g.).

Section 2 describes our straightforward computational scheme for the simulation of the atoms of a dense gas (listed in Ancillary Material, §A). Such a scheme is also at the core of more complicated schemes for more complicated molecular simulations (Koumoutsakos 2005; Horstemeyer 2009; Wagner et al. 2011, e.g.). In many scenarios it is easiest to write a microscale simulator with spatially periodic boundaries: for molecular dynamics some relevant comments by other researchers include "periodic boundaries have been used" (Evans and Hoover 1986, p.248); "In general, one prefers periodic boundary conditions" (Koplok and Banavar 1995, p.260); and "To circumvent this problem, ... a periodic system may be assumed" (Koumoutsakos 2005, p.477). One aim of the equation-free approach is to use whatever simulator has been provided, and adapt it to macroscale simulations. Hence the important new challenge we address is to use a triply-periodic atomic simulation code as 'a given' for the computed patches in an equation-free scheme.

In order to research realistic problems in the future we expect to implement the methodology within one of the established powerful molecular dynamics simulators such as LAMMPS (Plimpton et al. 2016). However, here we focus attention to establishing a proof-of-principle and the fundamental effectiveness of the scheme: for that purpose, our straightforward atomistic code is sufficient.

Our innovation is to couple relatively small, triply-periodic, atomistic patches over unsimulated space: the coupling mechanism has to differ from that used in patch schemes to date (Liu et al. 2015, e.g.). Section 3 describes a way to implement a proportional controller (Bechhoefer 2005, e.g.) in two so-called 'action regions' that surround the 'core' of each patch. The average kinetic energy in a core estimates the local temperature in a patch. Then interpolating such core temperatures over the unsimulated space estimates the macroscale temperature field. The applied control aims to appropriately drive the average kinetic energy in each action region to the corresponding macroscale temperature. Section 5 uses some modern dynamical systems theory (Aulbach and Wanner 2000) to prove that a line of such coupled periodic-patches has macroscale dynamics that emerge for a range of initial conditions and for a wide class of microscale systems. Then subsection 5.2 explicitly constructs the emergent dynamics for the scheme of controlled patches for a general advection-diffusion PDE. The construction establishes that there is a good control strength so that the emergent dynamics of the scheme reasonably approximates the correct macroscale advection-diffusion.

Section 4 confirms that the proposed controlled coupling of periodic patches is effective for atomistic simulations, and for a control roughly as predicted by the analysis.

For patch dynamics in space-time, a full implementation involves projective integration forward in time (Samaey, Roose, and Kevrekidis 2006; Givon, Kupferman, and Stuart 2004; Moller et al. 2005, e.g.). However, we leave projective integration of periodic patches to future research. Future research could also extend the analysis herein to establish the potential for high order accuracy, in multiple dimensions, analogous to what has been proven for patches with Dirichlet/Neumann/Robin boundaries (Roberts and Kevrekidis 2007; Roberts, MacKenzie, and Bunder 2014). Although this article focusses on the macroscale temperature diffusion emerging from an atomistic simulation, the equation-free patch scheme does usefully apply to wave systems (Cao and Roberts 2013; Cao and Roberts 2016) and so we expect that controlled periodic-patches should also be able to reasonably predict the emergent density-momentum waves of an atomistic simulation.

## 2 An isolated triply-periodic patch

As a first test of our novel methodology, our microscale, detailed, simulator is the molecular dynamics of a monatomic gas in 3D space. The simulator computes the motions of $N$ interacting atoms in a microscale patch of space-time, where for our purposes typically there are up to a few thousand atoms in a patch. For example, Figure 1 shows the apparently chaotic path in space of $N = 64$ atoms in a patch for one short-time simulation. We implemented a triply-periodic cubic domain where an atom crossing any face is re-injected into the cube across the opposing face. Our challenge is to develop methods that use such a 'given' spatially periodic microscale simulator to predict macroscale dynamics.

As evident in Figure 1, throughout we non-dimensionalise all quantities with respect to atomic scales so that, for example, the inter-atomic equilibrium distance is one, and the atomic mass is one.

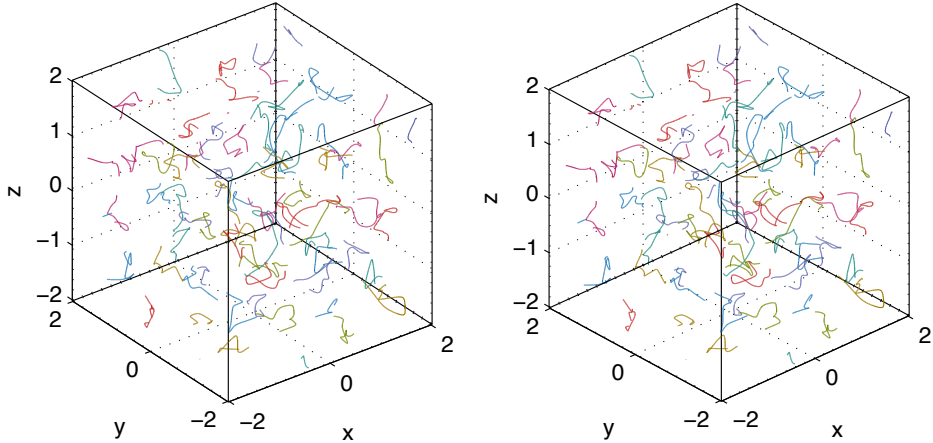

Figure 1: trajectories of 64 atoms, over a time $0 \leq t \leq 3$, in a triply-periodic, cubic, spatial domain, showing the beginnings of the complicated inter-atomic interactions. This stereo pair when viewed cross-eyed gives a 3D effect.

To describe the coded simulator (§A), let $\vec{x}_i(t)$ denote the position in space of the $i$th atom as a function of time $t$, and let $\vec{q}_i(t)$ denote the velocity of the $i$th atom. Then one set of ordinary differential equations (ODEs) for the system are (§A.3, lines 13–16)

$$\frac{d\vec{x}_i}{dt} = \vec{q}_i, \quad i = 1, \ldots, N. \tag{1}$$

The other set of ODEs for the system come from the inter-atomic interactions. For this monatomic gas we use the classic Lennard-Jones potential (Koplok and Banavar 1995, e.g.) for which the force between atoms separated by a distance $r$ is $F = 1/r^7 - 1/r^{13}$ (non-dimensionally). For atoms of non-dimensional mass one, Newton's 2nd law then gives the acceleration of each atom as (§A.3, lines 27–35)

$$\frac{d\vec{q}_i}{dt} = \sum_j \left( \frac{1}{r_{ij}^7} - \frac{1}{r_{ij}^{13}} \right) \frac{\vec{r}_{ij}}{r_{ij}}, \quad i = 1, \ldots, N, \tag{2}$$

where $\vec{r}_{ij}$ is the displacement vector from the $i$th atom to the $j$th, and distance $r_{ij} = |\vec{r}_{ij}|$. Because the patch is triply-periodic, the inter-atomic sum in (2) should be over all periodic images of the atoms. Computationally, in the sum we neglect atoms and their images further away than a patch half-width (§A.3, lines 19–24). Due to the $1/r^7$ decay of long range attraction, and with a typical patch of size $10 \times 10 \times 10$ atoms, the error in accounting for only these atoms/images is roughly $5^{-7} \approx 10^{-5}$—reasonably negligible. Upon checking, the mean momentum is conserved to machine precision which reflects the symmetry in the coding. Moreover, Figure 2 shows one example of the kinetic and potential energy, and their sum: we found the total energy is typically conserved to a relative error of about $10^{-5}$.

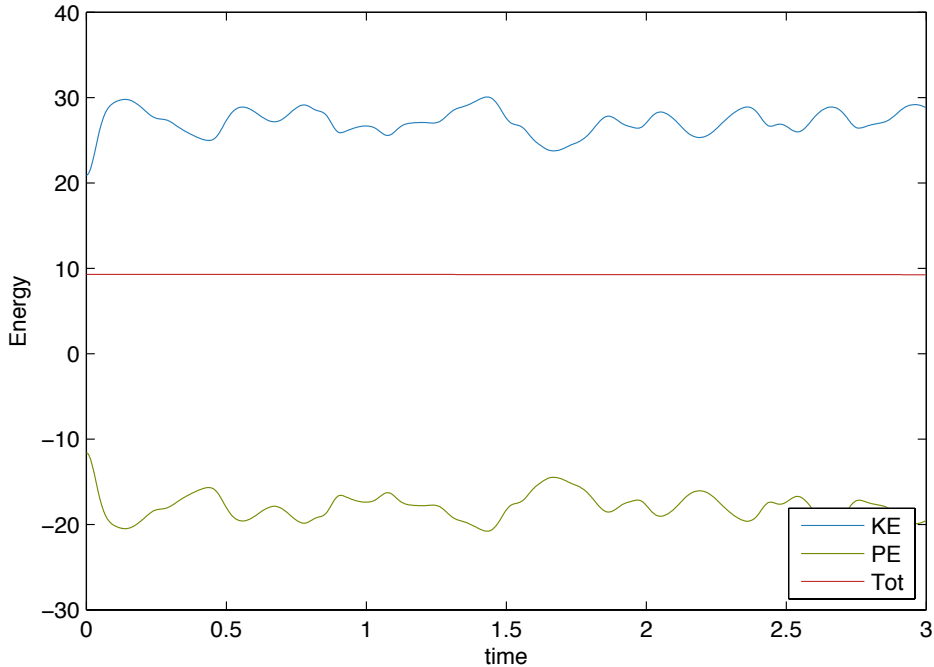

Figure 2: for the case of Figure 1, this plot of nondimensional Kinetic Energy (KE), Potential Energy (PE), and their sum (Tot) illustrates the conservation of energy by the coded simulation.

Most of the remaining code in the time derivative routine (§A.3) couples a patch to the macroscale surroundings—described by the next Section 3. Our code does not employ any fast multipole or cognate techniques (Cheng, Greengard, and Rokhlin 1999, e.g.) because preliminary exploration indicated that for only a few thousand atoms per patch the fast techniques are not effective. The microscale simulator is then to integrate the ODEs (1)–(2) in time. For simplicity in this proof-of-principle study, we use a generic MATLAB integration routine (§A.1, line 31) rather than any of the more accurate symplectic integrators that apply to this Hamiltonian system (Yoshida 1993; Hairer, Lubich, and Wanner 2003, e.g.). Such a generic integration routine would not handle the discontinuous reinjection of atoms that leave the cube: consequently, we allow the integrated atom positions to exit the cube smoothly, but map such atoms inside the cubic patch (§A.2) for computing inter-atomic forces and for plotting. Figure 1 plots one very short example simulation with $N = 64$ atoms.

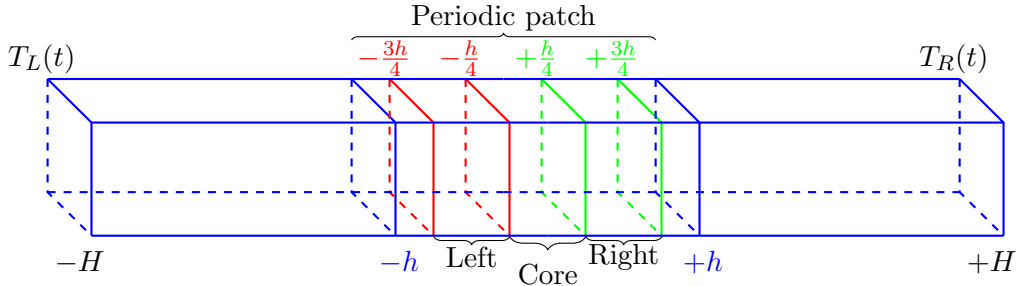

Figure 3: the simplest case is one triply-periodic patch of atomistic simulation, $-h < x < h$, coupled to distant sidewalls, at $x = \pm H$, of specified temperature. The patch's core region defines its local temperature, and a proportional controller is applied in the left and right action regions to generate a good macroscale prediction.

## 3 Couple patches with a proportional controller

The previous Section 2 described a microscale simulator for the isolated dynamics of atoms in a cubic domain. It is an example of the straightforward code that can be written for simulations which are triply-periodic in space. Our innovative challenge is to use the code, as if it were almost a 'black-box', to simulate over large space-time scales. In our equation-free methodology, such a large scale scape-time simulation is achieved by simulating in small, spatially distributed, patches and coupling them over the empty space between the patches (Samaey, Kevrekidis, and Roose 2005; Kevrekidis and Samaey 2009, e.g.).

In this first study of the use of periodic-patches, we only address the scenario of one large spatial dimension. The other space-time dimensions are assumed small. Further, as a proof-of-principle, this section addresses the specific case of one such periodic-patch coupled to 'distant' imposed boundary conditions (Section 5 analyses the case of multiple coupled patches in one large dimension). Figure 3 illustrates the scenario with one 'small' periodic-patch centred on $x = 0$ in a macroscale domain $-H < x < H$. This basic scenario allows us to focus on the key methodological innovation: namely how to couple a microscale periodic-patch to the surrounding macroscale environment.

For the purpose of validating our novel patch scheme, in this pilot study we suppose we want to predict macroscale heat transport by the atomistic simulation. Thus here we compare the scheme's predictions for the temperature field $T(x,t)$ with that of the continuum heat diffusion PDE

$$\frac{\partial T}{\partial t} = K\frac{\partial^2 T}{\partial x^2}, \quad \text{such that } T(-H,t) = T_L(t), \quad T(+H,t) = T_R(t). \qquad (3)$$

This section places at the origin a $(2h \times 2h \times 2h)$-periodic-patch of the atom simulation, as in Figure 3. This patch extends over $-h < x < h$ within the macroscale domain $-H < x < H$ with specified temperatures, $T_L(t)$ and $T_R(t)$, on the ends of the domain, $x = -H$ and $x = H$ respectively: un-simulated spaces are the comparatively large domains $h < |x| < H$. This article mostly uses the convention of lowercase letters denoting microscale quantities, such as $\vec{x}$, $\vec{q}$, and $h$, and uppercase letters denoting macroscale quantities, such as $H$ and $T$. Here the problem is assumed to be doubly-$2h$-periodic in the other two spatial dimensions; that is, the physical domain is long and thin (Figure 3). We divide the patch into four equal-sized regions (§A.3, lines 38–43):

- the core being $|x| < h/4$ that is used to define macroscale quantities of the patch such as the temperature $T(0,t)$;

- the left action region being $|x+h/2| < h/4$ that is, with the right action region, used to couple the patch to surround macroscale information such as the next patches or the environmental boundary values;

- the right action region being $|x - h/2| < h/4$; and

- a 'buffer' region for $|x| > 3h/4$ that caters for a smooth transition between the action regions.

Figure 4 shows an alternative schematic view of the patch: this view emphasises the microscale $2h$-periodicity in $x$ and the role of the 'buffer' region between the two action regions, on the 'opposite side' to the important core region.

There appears to be no need for buffer regions between the action regions and the core region (Bunder and Roberts 2012; Bunder, Roberts, and Kevrekidis 2017) (Figures 3 and 4).

We implement a proportional controller (Bechhoefer 2005, e.g.) to couple the patch to the surrounding macroscale variations (we leave exploration of proportional-integral and proportional-integral-derivative controllers to future research). The applied control is proportional to the differences in the action regions between the macroscale field and the microscale patch simulator. During the atomistic simulation we compute the nondimensional temperatures in the core and action regions as (§A.3, lines 45–48)

$$T_c = \operatorname*{mean}_{j \in \text{core}} \text{KE}_j\,, \quad T_l = \operatorname*{mean}_{j \in \text{left}} \text{KE}_j\,, \quad T_r = \operatorname*{mean}_{j \in \text{right}} \text{KE}_j\,, \qquad (4)$$

in terms of the non-dimensional kinetic energy of each atom, $\text{KE}_j = |\vec{q}_j|^2/2$ (the initial conditions and conservation of momentum in the algorithm ensure the mean velocity of the atoms is zero). With one patch centred at $x = 0$ coupled to boundaries at $x = \pm H$, the scheme's predicted macroscale field for the temperature is the parabolic interpolation through the three values $T_L$, $T_0$ and $T_R$,

$$T(x,t) = T_L \frac{x(x-H)}{2H^2} + T_0 \frac{H^2 - x^2}{H^2} + T_R \frac{x(x+H)}{2H^2}\,, \qquad (5)$$

where the small finite width of the core results in the central temperature $T_0 = T(0,t)$ being slightly different to the core mean $T_c$:

$$T_0 = \frac{T_c - (T_L + T_R)\frac{1}{6}(h/4H)^2}{1 - \frac{1}{3}(h/4H)^2}\,. \qquad (6)$$

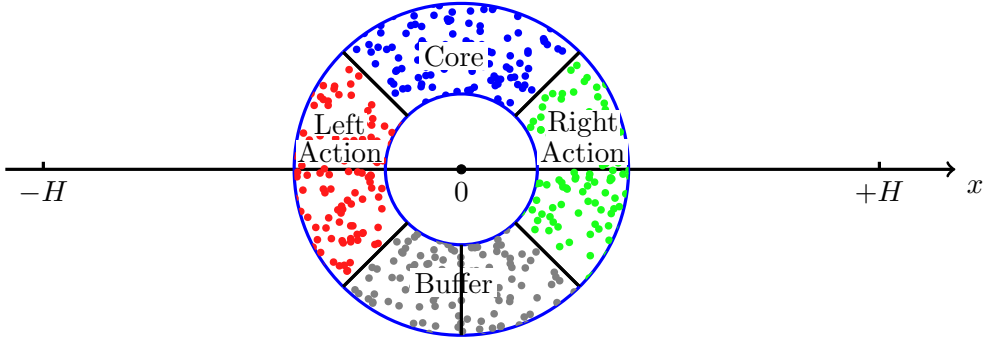

Figure 4: a schematic view of a microscale patch that emphasises the $2h$-periodicity in $x$ and indicating the need for the lower 'buffer' region allowing for a smooth transition between the action regions 'opposite' the core.

Then averaging the predicted macroscale field (5) over each action region gives that the macroscale interpolation predicts the action regions should have the average temperatures

$$T_{\text{int},\pm} = T_0 \pm (T_R - T_L)(h/4H) + \tfrac{13}{6}(T_R - 2T_0 + T_L)(h/4H)^2, \qquad (7)$$

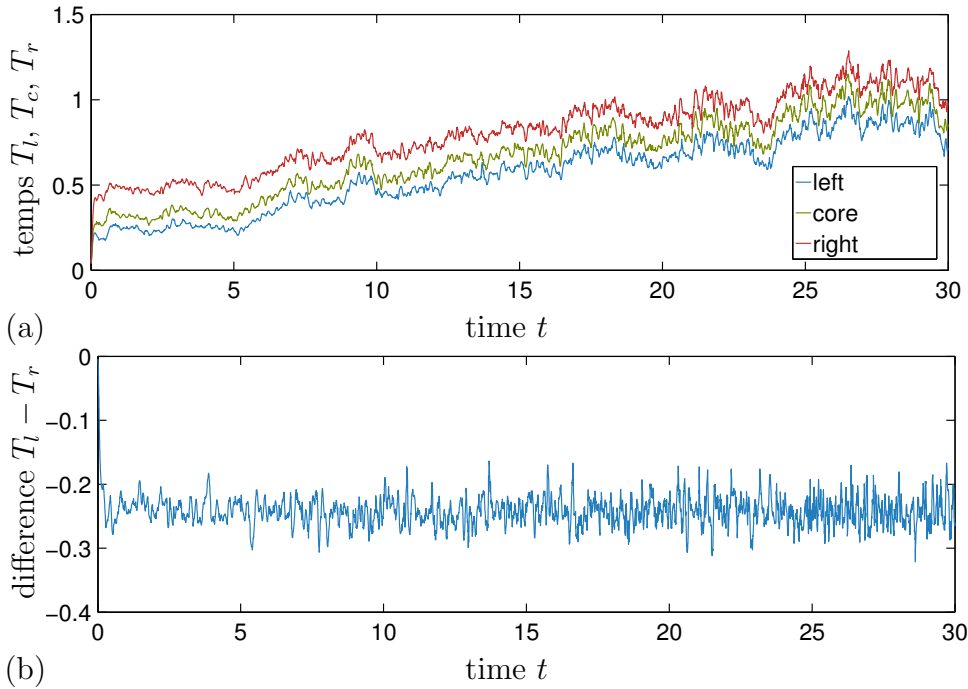

Figure 5: (a) temperatures over macroscale times in the sub-patch regions, and (b) the temperature difference $T_l - T_r$. The simulation is of 343 atoms in a patch of spatial periodicity $2h = 7$ and with control strength $\mu = 30$ to couple with macroscale boundary temperatures $T_R = 1.5$ and $T_L = 0.5$ at $x = \pm 7$.

where $\pm$ is for the right/left action region respectively. The applied control is proportional to the differences between these macroscale predictions (7) and the temperatures (4) from the patch simulation (§A.3, lines 50–55).

The controller accelerates or decelerates the atoms in the action region accordingly. That is, for each atom, Newton's second law (2) is modified by the control to (§A.3, lines 56–59)

$$\frac{d\vec{q}_j}{dt} = \cdots + \begin{cases} \frac{K\mu}{2h^2 T_r}(T_{\text{int},+} - T_r)\vec{q}_j, & j \in \text{right action}, \\ \frac{K\mu}{2h^2 T_l}(T_{\text{int},-} - T_l)\vec{q}_j, & j \in \text{left action}, \\ 0, & \text{otherwise}, \end{cases} \qquad (8)$$

where the dots denote the forces (2) from the interatomic Lennard-Jones potential, and nondimensional parameter $\mu$ is the strength of the control. Consequently, when the atoms in an action region are too cool, below $T_{\text{int},\pm}$, then the control accelerates the atoms to heat them up, and vice versa.

Figure 5 demonstrates the proportional controller is effective. The simulation is of 343 atoms in a patch of size $7 \times 7 \times 7$ coupled to boundaries at $x = \pm 7$ with specified temperatures $T_L = 0.5$ and $T_R = 1.5$. The control strength $\mu = 30$ in (8). Two time scales are apparent.

- On a microscopic time scale of $\Delta t < 1$ the control establishes that the temperatures in the action region differ according to the local gradient of macroscale temperature. In this realisation the temperature gradient is $1/(2H) = 1/14$, so that over the distance $2/7$ between mid-action regions a temperature difference of 0.25 is maintained as shown by Figure 5(b). The atomistic fluctuations about this mean increase in time as the overall temperature increases (Figure 5(a)).

- Figure 5(a) shows that over a macroscopic time scale of $\Delta t \approx 30$ the core temperature evolves towards the correct mean of 1.0 —albeit with fluctuations arising from the stochastic nature of the atomistic dynamics. This macroscopic time-scale is the diffusion time for heat across the macroscopic length scale $2H = 14$. Indeed it is on this and longer macroscopic time scales that in future developments we would implement projective integration in time (Gear and Kevrekidis 2003; Kevrekidis and Samaey 2009, e.g.).

For completeness, we should perhaps also implement proportional controllers for the other macroscale variables of density, pressure and average velocity. However, here there is little to control in these variables as, through conservation in the code of atoms and overall momentum, in the patch the density is constant, and the average velocity zero from the initial conditions. Thus in this first study we just control the nontrivial dynamics of the temperature.

## 4  Analyse optimal control for a single patch

This section analyses heat diffusion PDEs that model the controlled patch atomistic simulation of the previous Section 3. We find that a good control has strength parameter of roughly $\mu \approx 30$ in this scenario.

We compare the patch atomistic simulation with the dynamics of the well-established continuum heat diffusion PDE (3) on the macroscale domain. With constant boundary temperatures, its equilibrium solution is the linear field $T = \frac{1}{2}(T_R + T_L) + \frac{x}{2H}(T_R - T_L)$. The dynamics superimposed on this equilibrium are an arbitrary linear combination of the modes $e^{-\lambda_n t} \sin[k_n(x + H)]$ for eigenvalue $\lambda_n = -Kk_n^2$ and wavenumbers $k_n = n\pi/(2H)$. With only one patch in the domain (Figure 3), the patch atomistic simulation can only approximate the gravest $n = 1$ mode $e^{-\lambda_1 t} \cos(k_1 x)$ with wavenumber $k_1 =$ $\pi/(2H)$ and rate $\lambda_1 = -Kk_1^2 = -K\pi^2/(4H^2)$. We aim for the predictions (§3) of the single patch scheme to match this gravest mode.

The atomistic simulation within a patch will also be reasonably well modelled by the continuum heat diffusion PDE albeit with significant microscale fluctuations. In the microscale patch scheme here, Figure 3, the problem is homogeneous in the cross-sectional variables $y$ and $z$ so we only explore the $xt$-structure, $T(x,t)$. The continuum PDE for the controlled patch is then

$$\frac{\partial T}{\partial t} = K\frac{\partial^2 T}{\partial x^2} + \frac{K\mu}{h^2}g(x,T), \quad 2h\text{-periodic in } x, \tag{9}$$

where the control shape, piecewise constant, is

$$g(x,T) = \begin{cases} T_{\text{int},+} - \frac{2}{h}\int_{h/4}^{3h/4} T\,dx, & \frac{h}{4} < x < \frac{3h}{4}, \\ T_{\text{int},-} - \frac{2}{h}\int_{-3h/4}^{-h/4} T\,dx, & -\frac{3h}{4} < x < -\frac{h}{4}, \\ 0, & \text{otherwise}. \end{cases} \tag{10}$$

The temperatures $T_{\text{int},\pm}$, given by (7), come from the macroscale coupling which here is the parabolic interpolation (5)–(6) via the core average $T_c = \frac{2}{h}\int_{-h/4}^{h/4} T\,dx$ to $T_R$ and $T_L$ at boundaries $x = \pm H$.

In the case of constant boundary temperatures, this controlled patch problem has equilibrium, with $T_c = T_0 = \frac{1}{2}(T_L + T_R)$, of

$$T = \begin{cases} T_0 + \frac{\mu\Delta T}{4}x/h, & |x| < \frac{h}{4}, \\ T_0 + \frac{\mu\Delta T}{16}\left[\frac{3}{2} - 2(2x/h - 1)^2\right]\operatorname{sgn} x, & \frac{h}{4} < |x| < \frac{3h}{4}, \\ T_0 + \frac{\mu\Delta T}{4}(\operatorname{sgn} x - x/h), & |x| > \frac{3h}{4}, \end{cases} \tag{11}$$

where $\Delta T = (T_R - T_L)h/(4H)/[1 + \mu/12]$. A first important check is that this equilibrium has the correct mean core temperature of $T_0 = (T_R + T_L)/2$.

To understand the dynamics of the patch scheme (9)–(10) we characterise the general dynamics, relative to the equilibrium, in terms of spatial eigenfunctions. We express the dynamic analysis in terms of the half-width of the patch, $h$, and the ratio $r = h/H$ which is the fraction of macroscale space on which the microscale patch simulation is executed/solved. There are two classes of eigenfunctions: symmetric and antisymmetric within the patch.

- Symmetric eigenfunctions corresponding to temporal decay $e^{-Kk^2 t}$, $k > 0$, are of the form

$$T = \begin{cases} A\cos kx, & |x| < h/4, \\ C + D\cos[k(x \mp \frac{h}{2})] \pm E\sin[k(x \mp \frac{h}{2})], & \pm h/4 \lessgtr x \lessgtr \pm 3h/4, \\ B\cos[k(x \mp h)], & x \gtrless \pm 3h/4. \end{cases} \tag{12}$$

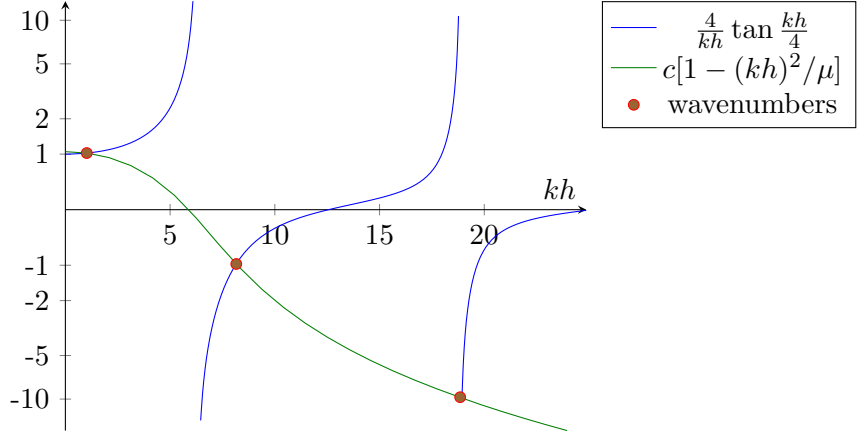

Figure 6: (blue) $(4/kh)\tan(kh/4)$ has an infinite number of vertical asymptotes at $kh = 2n\pi$ for odd $n$; (green) the parabola $c[1-(kh)^2/\mu]$ intersects it infinitely often (here for $\mu = 34.67$ and $r = 0.64$). The vertical axis is nonlinearly scaled.

Substitute these eigenfunctions into the governing equations (6), (7), (9) and (10). Then straightforward algebra, detailed by Alotaibi (2017, §2.4.2), leads to requiring the characteristic equation

$$\cos\frac{kh}{2}\sin\frac{kh}{4}\left[\left(\frac{1-7r^2/48}{1-r^2/48}\right)\frac{4}{kh}\sin\frac{kh}{4}+\left(\frac{1}{\mu}k^2h^2-1\right)\cos\frac{kh}{4}\right]=0\,. \tag{13}$$

The first factor being zero gives $kh/2 = n\pi/2$ for odd $n$; that is, wavenumber $k = n\pi/h$ for odd $n$. The second factor being zero gives $kh/4 = n\pi$ for integer $n$; that is, wavenumber $k = 4n\pi/h$ for integer $n$. The third factor being zero rearranges to

$$\frac{4}{kh}\tan\frac{kh}{4} = \left(\frac{1-r^2/48}{1-7r^2/48}\right)\left[1-\frac{1}{\mu}(kh)^2\right]. \tag{14}$$

Figure 6 plots the two sides of this equation illustrating that there are an infinite number of wavenumbers $k$, only one of which is small. It is this small wavenumber mode that is of macroscale interest.

- Antisymmetric eigenfunctions are of the form

$$T = \begin{cases} A\sin kx\,, & |x| < h/4\,, \\ \pm C + D\sin[k(x\mp\tfrac{h}{2})]\pm E\cos[k(x\mp\tfrac{h}{2})], & \pm h/4 \lessgtr x \lessgtr \pm 3h/4\,, \\ B\sin[k(x\mp h)], & x \gtrless \pm 3h/4\,. \end{cases}$$

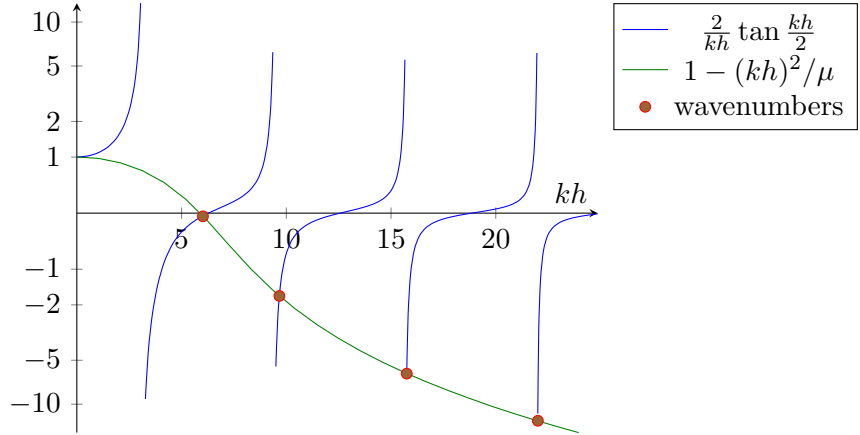

Figure 7: (blue) $(2/kh)\tan(kh/2)$ has an infinite number of vertical asymptotes at $kh = n\pi$ for odd $n$; (green) the parabola $1 - (kh)^2/\mu$ intersects it infinitely often (here for $\mu = 34.67$ and $r = 0.64$). The vertical axis is nonlinearly scaled.

Substituting into the governing equations and straightforward algebra, detailed by Alotaibi (2017, §2.4.3), leads to requiring the characteristic equation

$$\sin\frac{kh}{2}\left[\frac{2}{kh}\sin\frac{kh}{2} + \left(\frac{1}{\mu}k^2h^2 - 1\right)\cos\frac{kh}{2}\right] = 0. \qquad (15)$$

The first factor being zero gives $kh/2 = n\pi$ for integer $n$; that is, wavenumber $k = 2n\pi/h$ for integer $n$. The second factor being zero rearranges to

$$\frac{2}{kh}\tan\frac{kh}{2} = 1 - \frac{1}{\mu}(kh)^2.$$

Figure 7 plots the two sides of this equation illustrating that there are an infinite number of possible wavenumbers $k > 0$ (the algebraic formula degenerates at $k = 0$ so that apparent intersection is not a possible wavenumber).

With one exception, all of these possible wavenumbers $k \propto 1/h$, and the corresponding decay rate $\propto K/h^2$ characterises microscale modes within the patch. All these microscale modes cause the spatial structure internal in the patch to rapidly approach a quasi-equilibrium: in a time proportional to a cross-patch diffusion time of $h^2/K$.

The one exceptional mode corresponds to the small wavenumber $k$ shown in Figure 6. A straightforward small $kh$ approximation to the left-hand side

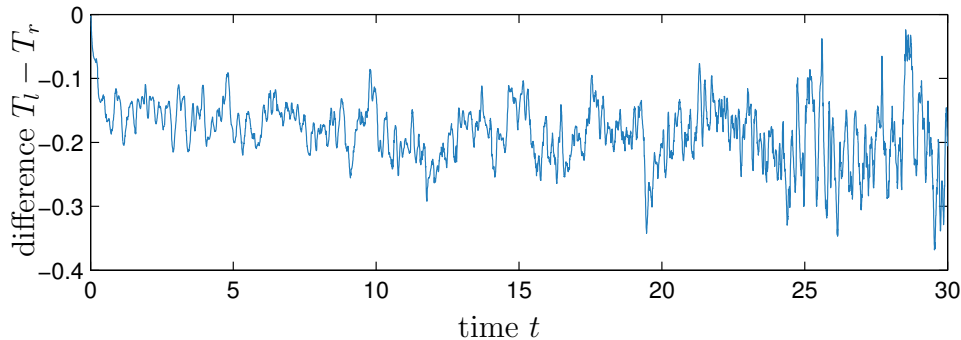

Figure 8: the temperature difference $T_l - T_r$. The simulation is of 343 atoms in a patch of spatial periodicity $2h = 7$ and with control strength $\mu = 3$ to couple with macroscale boundary temperatures $T_R = 1.5$ and $T_L = 0.5$ at $x = \pm 7$.

of (14) (Alotaibi 2017, §2.5) leads to

$$\left(\frac{kh}{4}\right)^2 \approx \frac{3r^2/8}{(1 - 7r^2/48) + (1 - r^2/48)48/\mu}.$$

Consequently $kh \propto r$. But the neglected terms in this approximation are $\mathcal{O}\big((kh)^4\big)$ which are equivalently $\mathcal{O}\big(r^4\big)$ and so for consistency the $\mathcal{O}\big(r^2\big)$ terms in the denominator should be neglected to lead to

$$kH = kh/r \approx \sqrt{\frac{6}{1 + 48/\mu}}. \qquad (16)$$

This wavenumber $\propto 1/H$ is characteristic of a macroscale mode, and the corresponding decay rate $\propto K/H^2$ also characterises a single macroscale mode across the domain. For the controlled patch simulation to best predict the correct macroscale dynamics of this mode, we need the wavenumber (16) to match the gravest wavenumber of the heat PDE (3) on the macroscale domain: the start of this section found it to be $k_1 H = \pi/2$. Thus best predictions, via a little algebra, suggest

$$\frac{1}{\mu} \approx \frac{24/\pi^2 - 1}{48} = 0.02983, \quad \text{that is,} \quad \mu \approx 33.53. \qquad (17)$$

for the control strength.

In an application, the result (17) requires an estimate of the diffusivity $K$. Simulations can give a rough estimate of diffusivity $K$ for any microscale system that is diffusive-like on the microscale. Figure 8 is for the same scenario as Figure 5, but with a weaker control strength of $\mu = 3$. Figure 8 suggests that the initial microscale transients decay on a time scale of roughly one (and then evolve slowly on a macroscale time). That is, the leading antisymmetric microscale wavenumber $k_3$ is such that $Kk_3^2 \approx 1$, that is, $K \approx 1/k_3^2$. Figure 7 illustrates that $\pi < k_3 h < 2\pi$ depending upon the control via $1/\mu$. Assuming $k_3 \ell \approx 3\pi/2$—mid-range should be good enough for rough estimates—then the microscale diffusion constant $K \approx 4h^2/(9\pi^2)$. The simulations of Figures 5 and 8 were in a patch of nondimensional length $2h = 7$, hence here the nondimensional diffusion constant is $K \approx 0.5$, roughly. This resultant control is roughly what we found convenient for generating Figures 5 and 8. The importance of this paragraph is that it illustrates how a little analysis and a few simulations can determine a reasonably good control of the microscale periodic patch.

## 5 Couple multiple periodic patches in general

The aim of this section is to establish a basis for theoretical support of the controlled patch/gap-tooth scheme for general molecular/particle/agent based simulations. Consider the class of systems whose mesoscale dynamics are modelled by a field $u(x,t)$ governed by the stochastic reaction-advection-diffusion PDE

$$\frac{\partial u}{\partial t} = \mathcal{L}u + \alpha f(u, u_x) + \sigma g(u, u_x)\dot{W}(x,t). \tag{18}$$

This stochastic PDE (18) (SPDE) is termed "mesoscale" because it codifies the in-principle closure of some atomistic simulation on length-time scales large enough for a reasonable closure to exist (even when unknown), but smaller than the macroscale of interest, and small enough so that chaotic atomistic fluctuations may appear as noise $\dot{W}$ of strength $\sigma$. The functions $f$ and $g$ are some functions for the closure of this mesoscale model. The closure functions $f$ and $g$ need not be known for our theoretical support to apply, but we do assume that $f$ and $g$ are smooth enough. Here we focus attention to systems where the linear operator satisfies three properties: $\mathcal{L}(\text{constant}) = 0$; $\mathcal{L}$ is otherwise dissipative on the mesoscale–macroscale, such as $\mathcal{L} \approx K\partial_{xx}$; but where $\mathcal{L}$ is bounded due to 'cut-offs' on the microscale. This section invokes stochastic centre manifold theory (initiated by Boxler (1989)) to provide novel support for the contention that controlled patches of the SPDE (18) successfully simulates its macroscale dynamics. Consequently, for a wide range of microscale simulators the same control of periodic patches will economically predict the macroscale dynamics.

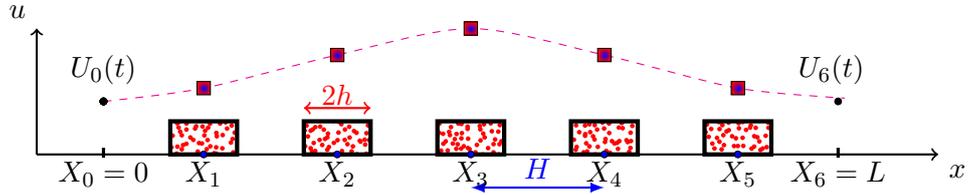

Figure 9: schematic illustration of five microscale patches, $2h$-periodic, centred at $X_j = jH$, inside a macroscale domain $[0, L]$. The field $u(x,t)$, satisfying the SPDE (18), has some macroscale boundary conditions such as the Dirichlet conditions that $u(0,t) = U_0(t)$ and $u(L,t) = U_6(t)$.

This section does not precisely prescribe all restrictions on the form of the SPDE (18) because the aim of this section is to scope the domain and feasibility of the controlled periodic-patch scheme. The issue of whether specific systems satisfy the precise requirements for rigorous support is left for future researchers to certify for the systems of their interest.

## 5.1 Macroscale existence and emergence theory

Figure 9 illustrates a typical scenario. Suppose we are interested in the dynamics of the field $u$ on a relatively large spatial domain $\mathbb{X} = [0, L]$. The patch scheme distributes $M$ patches, equi-spaced with macroscale spacing $H$, and centred at $X_j = jH$. The microscale patches are $2h$-periodic with relatively small half-size $h$: the scale ratio $r = h/H$ would be small for efficient simulation. Let the $2h$-periodic field in the $j$th patch be denoted by $u_j(x,t)$ for $|x - X_j| < h$ —distinct from the field $u(x,t)$ satisfying the SPDE (18) over all $\mathbb{X}$.

The prime quantity of interest on the macroscale the simulation is a measure of the field in each patch. Let the overall field in the $j$th patch be measured by the average over the $j$th core region $|x - X_j| < h/4$:

$$U_j(t) := \frac{2}{h} \int_{X_j - h/4}^{X_j + h/4} u_j(x,t)\,dx\,, \quad j = 1, \ldots, M\,. \tag{19}$$

With $M$ patches, specified Dirichlet boundary values at $x = 0, L$ may be referred to by synonyms $U_0(t)$ and $U_{M+1}(t)$. Let collective core averages by denoted by the vector $\vec{U} = (U_1, U_2, \ldots, U_M)$.

Each periodic patch is controlled by some coupling with neighbouring patches. Suppose the simulation in each patch has an imposed proportional controller so the resultant effective system in the patches, $|x - X_j| < h$, modifies the effective SPDE (18) to

$$\begin{aligned}\frac{\partial u_j}{\partial t} &= \mathcal{L}u_j + \alpha f(u_j, u_{jx}) + \sigma g(u_j, u_{jx})\dot{W}_j(x,t) \\ &\quad + \mu\left\{\left[I_j^+(\vec{U},\gamma) - u_j^+\right]\chi_j^+(x) + \left[I_j^-(\vec{U},\gamma) - u_j^-\right]\chi_j^-(x)\right\}. \end{aligned} \quad (20)$$

The strength of the control is parametrised by $\mu$, and is applied in the left and right action regions of each patch as coded by the indicator functions

$$\chi_j^\pm(x) := \begin{cases} 1, & h/4 < \pm(x - X_j) < 3h/4\,, \\ 0, & \text{otherwise.} \end{cases}$$

The control depends upon the difference between an interpolation, $I_j^\pm(\vec{U},\gamma)$, of the core averages in neighbouring patches and the local averages of the patch field in the left and right action regions,

$$u_j^\pm(t) := \frac{2}{h}\int_{X_j - h}^{X_j + h} u_j(x,t)\chi_j^\pm(x)\,dx\,.$$

Because Lagrange interpolation is known to produce systematic consistency in simpler scenarios (Roberts and Kevrekidis 2007; Roberts, MacKenzie, and Bunder 2014), to-date we have coupled periodic patches via the Lagrange interpolation of $\vec{U}$ to estimate $U$ at the mid-action points $x = X_j \pm \frac{1}{2}h$:

$$I_j^\pm := \left\{1 + \gamma\left[\pm\tfrac{1}{2}r\bar{\mu}\delta + \tfrac{1}{8}r^2\delta^2\right] + \gamma^2\left[\mp r(\tfrac{1}{4} - \tfrac{1}{8}r)\bar{\mu}\delta^3 - r(\tfrac{1}{8} - \tfrac{1}{16}r)\delta^4\right] + \cdots\right\}U_j\,, \quad (21)$$

expressed in terms of centred difference and mean operators $\delta U_j = U_{j+1/2} - U_{j-1/2}$ and $\bar{\mu}U_j = (U_{j+1/2} + U_{j-1/2})/2$. The coupling parameter $\gamma$ is an artificial homotopy parameter that empowers us to connect a theoretical base at $\gamma = 0$ to the fully coupled case $\gamma = 1$ which is of interest. That is, this embeds the relevant physical problem, at $\gamma = 1$, into a family of problems parametrised by $\gamma$. For the theory of this section we require three things of the dependency of the interpolation operator $I_j^\pm$ on parameter $\gamma$: $I_j^\pm(\vec{U},\gamma)$ be smooth; $I_j^\pm(\vec{U},0) = U_j$; and $I_j^\pm(\vec{U},1)$ corresponds to the coded coupling of the computational patch scheme. Further, although not necessary, it is convenient to express $I_j^\pm(\vec{U},\gamma)$ as a polynomial in $\gamma$, as in (21), such that a truncation which neglects terms $\mathcal{O}(\gamma^{p+1})$ invokes an $I_j^\pm$ which depends upon only $U_{j-p}, \ldots, U_{j+p}$. This property then empowers convenient comparison with classic finite differences/elements.

Theoretical support is based upon an equilibrium of the controlled patch system (20). For parameters $\alpha = \sigma = \gamma = 0$ the system (20) is linear with $u_j(x,t) = $ constant being equilibria. Being independently constant in each of the $M$ patches, these equilibria form an $M$-D subspace $\mathbb{E}_0$ (a subspace of the space of $(\alpha, \sigma, \gamma, u_1(x), u_2(x), \ldots, u_M(x))$). By the definition (19) of the core averages, we write these equilibria as the patch field $u_j(x,t) = U_j$ such that $dU_j/dt = 0$ for $j = 1, \ldots, M$.

This subspace commonly and usefully persists upon perturbation. Specifically, we are interested in perturbations by nonlinearity (non-zero $\alpha$), by stochasticity (non-zero $\sigma$), and most importantly by coupling with neighbouring patches (non-zero $\gamma$). Consider the dynamics of the controlled patches, SPDE (20), linearised about each of the equilibria in $\mathbb{E}_0$, namely, for patch field $u_j$ being $2h$-periodic,

$$\frac{\partial u_j}{\partial t} = \mathcal{L} u_j + \mu \left\{ \left[ U_j - u_j^+ \right] \chi_j^+(x) + \left[ U_j - u_j^- \right] \chi_j^-(x) \right\}. \qquad (22)$$

As for diffusion, $\mathcal{L} = K \partial_{xx}$, we assume that the operator on the right-hand side of (22) has spectrum $0 = \lambda_1 > -\beta > \Re\lambda_2 \geq \Re\lambda_3 \geq \cdots$, where the size, $\beta$, of the spectral gap depends upon the size of the patches. In the linearisation (22), each of the $M$ patches are isolated from each other and consequently all of these eigenvalues have multiplicity $M$ (at least). For example, in the case of mesoscale diffusion, $\mathcal{L} = K \partial_{xx}$, the analysis of Section 4 applies except with the replacement of the factor $(1 - 7r^2/48)/(1 - r^2/48)$ by simply one: consequently Section 4 establishes that here $\lambda_1 = 0$, $\lambda_2 = -K\pi^2/h^2$, and so on, and consequently we may choose bound $\beta = 9K/h^2$, say. There are three extra zero eigenvalues, one each for the parameters $\alpha$, $\sigma$ and $\gamma$, corresponding to the extended state space formed by adjoining $d\alpha/dt = d\sigma/dt = d\gamma/dt = 0$. In this scenario, and subject to some technical caveats on the nature of the functions $f(u, u_x)$ and $g(u, u_x)\dot{W}(x, t)$, the marvellous theory of Aulbach and Wanner (2000) applies to establish that there exists an $M + 3$ dimensional slow manifold to the controlled patch mesoscale system (20) (or $mM + 3$ dimensional when the zero eigenvalue of (22) has multiplicity $m$ within each patch). The slow manifold exists globally in $\vec{U}$, and in a finite domain about zero in the parameters $(\alpha, \sigma, \gamma)$.[1]

Further, the theory establishes that the slow manifold emerges exponentially quickly from all initial conditions in some finite domain about the slow manifold (Aulbach and Wanner 2000, §4). More specifically, the exponential transients are $\mathcal{O}(e^{-\beta t})$ as $t \to \infty$. These transients are all microscale sub-patch modes as they correspond to the dissipative modes within isolated patches. For example, in the diffusive case the lower bound $\beta = 9K/h^2$ on the decay rate corresponds to an intra-patch diffusion time which is very small for the useful case of small patches. Just as for the single patch case of Section 4, this leaves the slow manifold and the evolution thereon to be the emergent macroscale dynamics of the controlled periodic-patch scheme.

---

[1] Although mesoscale diffusion $K\partial_{xx}$ is strictly an unbounded operator, such a diffusive operator is only a long wavelength approximation to the microscale dynamics which usually has a bounded equivalent operator (e.g., difference operators on a lattice). Hence the simulation $\mathcal{L}$ is usually bounded as required by Aulbach and Wanner (2000). The global existence explicitly proved by Aulbach and Wanner (2000) also requires perturbations to be Lipschitz and bounded, whereas here the inter-patch coupling as posited is unbounded: consequently, here we establish theoretical support in a finite local domain instead of the global domain.

## 5.2 Constructing a slow manifold predicts the macroscale

The previous subsection 5.1 establishes that the controlled periodic-patch scheme has emergent macroscale dynamics. To confirm, or otherwise, that the patch scheme's macroscale dynamics match that of the simulation, we have to construct the slow manifold and its evolution and compare with the mesoscale SPDE (18). Because the details of the slow manifold construction are closely linked to the details of the underlying system, here we specifically address the class of systems whose mesoscale has the diffusive $\mathcal{L} = K\partial_{xx}$. Further, noise from microscale chaos enormously complicates the construction (Roberts 2006; Bunder and Roberts 2017, e.g.), so here we further restrict analysis to cases with $\sigma = 0$.

A good robust method to construct slow manifolds is an iteration based upon the residuals of the governing equations (Roberts 1997; Roberts 2015, Part V): successive approximations to the slow manifold are corrected until the residuals are zero to within a pre-specified order of error. The algebraic machinations of the construction are involved and of little interest so are not recorded here. Instead we comment on two aspects:

- computer algebra readily constructs a slow manifold model of the macroscopic dynamics of the controlled periodic-patch scheme (Ancillary Material, §B) for a variety of systems in the class (18);

- as an example, we show that the controlled periodic-patch scheme applied to an advection-diffusion reasonably approximates the correct macroscale advection-diffusion.

For example, consider a microscale simulation whose mesoscale is the linear advection-diffusion PDE

$$u_t = Ku_{xx} - \alpha u_x, \tag{23}$$

with macroscale boundary conditions of $(MH)$-periodicity so that by symmetry each controlled patch is identical. Computer algebra, detailed in the

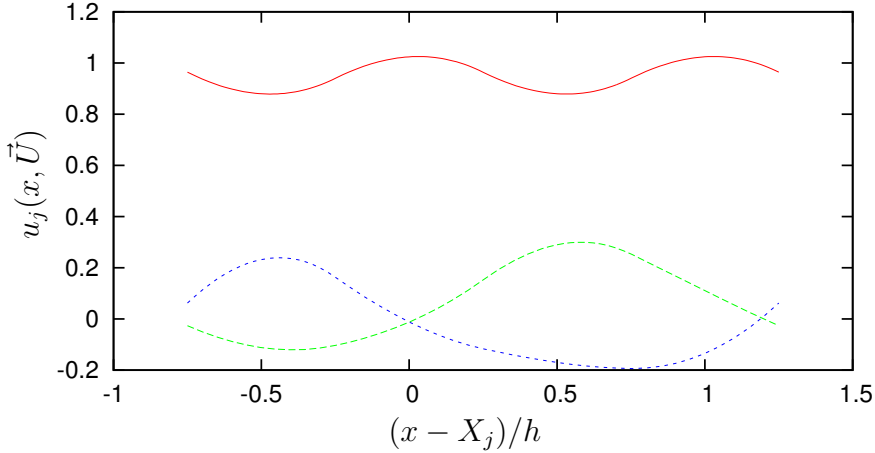

Figure 10: subpatch periodic field $u_j(\xi, \vec{U})$ for amplitudes $\vec{U} = \vec{0}$ except for: red-solid, $U_j = 1$; green-dashed, $U_{j+1} = 1$; blue-dotted, $U_{j-1} = 1$. For illustrative purposes this case is for $r = h = H = 1$, $\mathcal{L} = \partial_{xx}$, advection $\alpha = 1$, and control $\mu = 30$. Consequently the patch core is $|x - X_j|/h < \frac{1}{4}$, and action regions are $\frac{1}{4} < |x - X_j|/h < \frac{3}{4}$.

Ancillary Material §B, constructs the slow manifold of the corresponding controlled periodic-patch system (20). As the simplest example, with no advection, $\alpha = 0$, three iterations finds that each sub-patch field is piecewise parabolic: in terms of $\xi = (x - X_j)/h$ the $j$th patch field

$$u_j = U_j + \gamma \left[ U_j \frac{r^2 - 48\xi^2}{16(1 + 48/\mu)} + \sum_{\pm} U_{j\pm 1} \frac{-r^2 \pm 24\xi + 48\xi^2 + 12(-r^2 \pm 96\xi + 48\xi^2)/\mu}{32(1 + 48/\mu)(1 + 12/\mu)} \right] + \mathcal{O}(\gamma^2),$$

in the core $|\xi| < \frac{1}{4}$, and somewhat longer expressions in the other regions. Figure 10 plots leading order shapes of the fields in a patch across all four regions in the case of advection $\alpha = 1$: the general field in the patch is the linear combination of these curves with coefficients of $U_j$ and $U_{j\pm 1}$. This example is just a specific case of the reasonably general construction of the macroscale slow manifold. The small asymmetry in Figure 10 is due to the advection at velocity $\alpha = 1$ in this case. The downward curvature of the red-solid curve in the core region, $|\xi| < \frac{1}{4}$, and the back/buffer region, is due to the out-of-equilibrium macroscopic decay in the $j$th patch when surrounded by patches with zero core average, $U_{j\pm 1} = 0$. However, the interest on the macroscale is not the detailed sub-patch fields, but the long-term evolution of the macroscale amplitudes $\vec{U}$.

The macroscale amplitudes evolve according to the dynamics on the slow manifold. Simultaneous with constructing the shape of the slow manifold, the computer algebra (Ancillary Material §B) also finds the evolution of $\vec{U}$ on the slow manifold as a system of coupled ODEs. For the general advection-diffusion mesoscale PDE (23) the evolution on the slow manifold is

$$\frac{dU_j}{dt} = -\gamma\alpha\frac{15(U_{j+1} - U_{j-1})}{8H(1 + 48/\mu)(1 + 12/\mu)} + \gamma K\frac{3(U_{j-1} - 2U_j + U_{j+1})}{H^2(1 + 48/\mu)} + \mathcal{O}(\gamma^2, \alpha^2).$$

Evaluating at full coupling $\gamma = 1$, this slow manifold analysis then predicts that the emergent macroscale dynamics of the controlled periodic-patches of advection-diffusion (23) is

$$\frac{dU_j}{dt} \approx -\alpha\frac{15(U_{j+1} - U_{j-1})}{8H(1 + 48/\mu)(1 + 12/\mu)} + K\frac{3(U_{j-1} - 2U_j + U_{j+1})}{H^2(1 + 48/\mu)}. \quad (24)$$

Via such a model, the theory of this section empowers us to predict what the controlled periodic-patch scheme would do for a range of mesoscale systems.

In this particular example of general advection-diffusion we use the predicted model (24) to now identify a good control strength in the periodic-patch scheme. The macroscale discrete system (24) has equivalent PDE

$$U_t = -\alpha\frac{15}{4(1 + 48/\mu)(1 + 12/\mu)}U_x + K\frac{3}{1 + 48/\mu}U_{xx} + \mathcal{O}(H^2).$$

Thus the effective advection of the periodic-patch scheme will equal that of the posed underlying mesoscale system (23) when the control strength is chosen such that

$$\frac{15}{4H(1 + 48/\mu)(1 + 12/\mu)} = 1, \quad \text{that is,} \quad \mu \approx \frac{120 + 24\sqrt{69}}{11} = 29.0326.$$

Similarly, the effective diffusion will equal that of the posed underlying mesoscale system (23) when

$$\frac{3}{1 + 48/\mu} \approx 1 \quad \text{that is,} \quad \mu \approx 24.$$

To this level of approximation there may not be one single best control strength. Nonetheless, the closeness of these two good control strengths is encouraging for applications of the controlled periodic-patch scheme. Systematic exploration to higher order, and with other controller schemes, would give more accurate predictions for such optimal control in a wider range of mesoscale systems.

A future extension to a wider range of mesoscale systems is important because the primary rationale for the equation-free methodology is that its main application is to systems for which we do not know a meso–macro-scale closure. Hence we need to develop control schemes useful for a wide range of systems such as the class (18).

# 6 Conclusion

The atomistic simulation described by Section 2 is just one important example of microscale simulators used widely in engineering and science. In particular we address the class of simulators that are given with periodic conditions on the microscale. The challenge is to create a wrapper around such microscale periodic simulators in order to effectively predict macroscale behaviour. Section 3 implemented a proportional controller applied to action regions in the patch to couple patches to its neighbours over unsimulated space.

Analysis of atomistic simulations based upon just one small patch, and the corresponding controlled diffusion PDE, Section 4, indicates there are good values for the control strength. Section 5 then creates a theoretical basis for certifying that the emergent behaviour of many coupled, controlled, periodic patches does indeed predict appropriate macroscale, system level, dynamics for a wide range of microscale simulators.

The construction of a slow manifold model for the equation-free patch scheme provides an innovative route to explore algebraically how best to couple such periodic patches. Here, subsection 5.2 established the algebraic analysis to be feasible. Further research could search for optimal core and action regions sizes, optimal weight functions for the averages in the regions, refining the interpolation that couples the patches, and to other controllers.

**Acknowledgements**  Parts of this work was supported by the Taif University in Saudi Arabia, and by grant DP150102385 from the Australian Research Council. We thank Yannis Kevrekidis for many useful discussions and support.

# A Ancillary material: code for 3D atomistic simulation

This ancillary material provides the numerical code for simulating a microscale patch of atoms coupled over macroscale empty space to boundary values of temperature. It is included to document and potentially reproduce the results.

## A.1 Main driver code

```matlab
% coded 3d simulation of position and velocity of
% interacting atoms.  Uses ode23 to do time integration.
% HA 2015--6--19
global ll mucontrol TL TR ii nfns nAux hh Khsq

nAtom=64  % number of atoms
tEnd=3 % end time of simulation
ll=nAtom.^(1/3);   % length of periodic patch (inter-atom eq is at one)
mucontrol=0 % zero is no control, circa 30 seems best
TL=1, TR=0.5 % macroscale boundary values of temp
hh=ll  % macroscale BCs applied at +/-hh
rng('shuffle'); seed=422 %100+floor(900*rand); %random realisation seed
fileroot=['ctrlpatch' num2str(seed) 'N' num2str(nAtom)]
nAux=12; % number of auxillary variables computed
Khsq=0.5/(ll/2)^2; % coefficient of control

rng(seed);
% distribute atoms, randomly up to one per box
ns=ceil(ll);
i=linspace(-0.5,0.5,2*ns+1);
[j,i,p]=meshgrid(ll*i(2:2:end));
[~,k]=sort(rand(size(i(:)))); k=k(1:nAtom);
xq=[i(k) j(k) p(k) zeros(nAtom,3)]';
% add smallish, mean-zero, random position and velocity
zz=rand(6,nAtom)-0.5; zz=zz-repmat(mean(zz,2),1,nAtom);
xq=xq+0.3*diag([1,1,1,2,2,2])*zz; xq=xq(:);
% for imposing triple periodicity in space
ii=1:6:6*nAtom; ii=[ii ii+1 ii+2]+nAux;

% simulation in time from given ICs
```

```matlab
31  nfns=0;
32  [ts,Txqs]=ode23(@ctrlhddtu3dode,[0 tEnd],[zeros(nAux,1);xq]);
33  % impose periodicity on the computed positions
34  xqs=Txqs(:,nAux+1:end);   ii=ii-nAux;
35  xqs(:,ii)=xqs(:,ii)-round(xqs(:,ii)/ll)*ll;
36  % auxillary quantities at middle of time steps
37  Tuvws=diff(Txqs(:,1:nAux))./repmat(diff(ts),1,nAux);
38  tsx=(ts(1:end-1)+ts(2:end))/2;
39  nFunctions=nfns
40
41  ctrlgraphs % draw graphical output
```

## A.2 Interpose periodicity on positions

This function avoids MATLAB's ode23 objecting to discontinuities as atoms move across edges of the periodic box.

```matlab
1  function dxq=ctrlhddtu3dode(t,xq)
2  % Computes time derivative of position and velocity of
3  % interacting particles for Matlab integrator.
4  % HA Jan 2015 -- 2016
5  global  ll ii nfns
6  nfns=nfns+1;
7  % impose triple periodicity on positions
8  xq(ii)=xq(ii)-round(xq(ii)/ll)*ll;
9  dxq=ctrlhddtu3d(xq,t);
10 end
```

## A.3 Time derivatives of position and velocity

```matlab
1  function dTxq=ctrlhddtu3d(Txq,t)
2  % Computes time derivative of position and velocity of
3  % interacting particles.  Force triply periodic in
4  % space, periodicity ll. For the moment ignore that this
5  % dynamical system is sympletic.  HA Jan 2014 -- 2016
6  global  ll mucontrol TL TR nAux hh Khsq
7  % unpack positions
8  xq=Txq(nAux+1:end);
9  x=xq(1:6:end);
10 y=xq(2:6:end);
11 z=xq(3:6:end);
```

```matlab
% d/dt position = velocity
dxq=nan(size(xq));
dxq(1:6:end)=xq(4:6:end);
dxq(2:6:end)=xq(5:6:end);
dxq(3:6:end)=xq(6:6:end);

% Assume triple images enough to capture all significant.
[xxt,xx,lls]=meshgrid(x,x,[-ll 0 ll]);
[yyt,yy,lls]=meshgrid(y,y,[-ll 0 ll]);
[zzt,zz,lls]=meshgrid(z,z,[-ll 0 ll]);
[ddx,px]=min((xx-xxt+lls).^2,[],3);
[ddy,py]=min((yy-yyt+lls).^2,[],3);
[ddz,pz]=min((zz-zzt+lls).^2,[],3);
ds=sqrt(ddx+ddy+ddz)+1e-8;
% forces as a function of distance (with extra / dist)
fs=ds.^(-7-1)-ds.^(-13-1);
fs=-min(100,-fs); % ad hoc limit on force
fx=(xxt(:,:,1)-xx(:,:,1)-(px-2).*ll).*fs;
fy=(yyt(:,:,1)-yy(:,:,1)-(py-2).*ll).*fs;
fz=(zzt(:,:,1)-zz(:,:,1)-(pz-2).*ll).*fs;
% d/dt velocities = sum of forces
dxq(4:6:end)=sum(fx,2);
dxq(5:6:end)=sum(fy,2);
dxq(6:6:end)=sum(fz,2);

% proportional controller of the patch temperature
halfcore=ll/8; % halfwidth of core and action regions
xaction=ll/4;  % action regions centred at quarter points
% Find which atoms are in each region
jl=(abs(x+xaction)<halfcore);
jc=(abs(x        )<halfcore);
jr=(abs(x-xaction)<halfcore);
% unpack velocities and KEs, to find regional temperatures
u=xq(4:6:end); v=xq(5:6:end); w=xq(6:6:end); ke=(u.^2+v.^2+w.^2)/2;
Tl=mean(ke(jl));
Tc=mean(ke(jc));
Tr=mean(ke(jr));
% Control towards specified environmental temperature
rlh=ll/8/hh;
T0=(Tc-(TR+TL)*rlh^2/6)/(1-rlh^2/3);
Tintl=T0-(TR-TL)*rlh+(TR-2*T0+TL)*rlh^2*13/6;
```

```
53  Tintr=T0+(TR-TL)*rlh+(TR-2*T0+TL)*rlh^2*13/6;
54  ratel=mucontrol*(Tintl-Tl)*Khsq/Tl/2;
55  rater=mucontrol*(Tintr-Tr)*Khsq/Tr/2;
56  % de/accelerate in dirn of velocity, propto control
57  dxq(4:6:end)=dxq(4:6:end)+(ratel*jl+rater*jr).*u;
58  dxq(5:6:end)=dxq(5:6:end)+(ratel*jl+rater*jr).*v;
59  dxq(6:6:end)=dxq(6:6:end)+(ratel*jl+rater*jr).*w;
60
61  % other auxillary quantities are mean velocity in regions
62  ul_av=mean(u(jl)); vl_av=mean(v(jl)); wl_av=mean(w(jl));
63  uc_av=mean(u(jc)); vc_av=mean(v(jc)); wc_av=mean(w(jc));
64  ur_av=mean(u(jr)); vr_av=mean(v(jr)); wr_av=mean(w(jr));
65  aux=[Tl;Tc;Tr;ul_av;vl_av;wl_av;uc_av;vc_av;wc_av;ur_av;vr_av;wr_av];
66
67  % return auxillary quantities and atomic time derivatives
68  dTxq=[aux;dxq];
69  end
```

# B   Ancillary material: computer algebra derives a slow manifold

The following computer algebra code constructs a macroscale slow manifold to the controlled periodic-patch mesoscale PDE (20) in the autonomous case $\sigma = 0$. We use the language Reduce[2] because it is both free and perhaps the fastest general purpose computer algebra system (Fateman 2003).

Make printing pretty.
```
1 on div; off allfac; on revpri;
2 factor gamma,mu,delta,hh,kk,alpha;
```
Scale space to $\xi = (x - X_j)/H$ where `hh` $= H$.
```
3 depend xi,x;
4 let df(xi,x)=>1/hh;
```
Parametrise solution in terms of amplitudes $U_j$ such that $\partial U_j/\partial t = g(\vec{U}, t)$. Let's define $U_j$ to be the average over an averaging region, then the interpolated values will just be the required averages in those regions: although there may be a glitch for physical boundary conditions.
```
5 operator uu; depend uu,t;
6 let df(uu(~k),t)=>sub(j=k,g);
```

---

[2]http://reduce-algebra.com/

Zeroth approximation field is constant in each element. Write the patch field in four quarters: the core $u_c(\xi,t)$ for $|\xi| < \frac{1}{4}r$; the left action region $u_l(\xi,t)$ for $-\frac{3}{4}r < \xi < \frac{1}{4}r$; the right action region $u_r(\xi,t)$ for $\frac{1}{4}r < \xi < \frac{3}{4}r$; and the 'back/buffer' field $u_b(\xi,t)$ for $\frac{3}{4}r < |\xi| \leq r$.

```
 7 uc:=ul:=ur:=ub:=uu(j);
 8 g:=0;
```

Define averages over various regions.

```
 9 operator meanc; linear meanc;
10 operator meanl; operator meanr;
11 qr:=r/4;
12 let { meanc(xi^~~p,xi)=>qr^p*(1+(-1)^p)/2/(p+1)
13      , meanc(1,xi)=>1
14      , meanr(~a,xi)=>meanc(sub(xi=+r/2+xi,a),xi)
15      , meanl(~a,xi)=>meanc(sub(xi=-r/2+xi,a),xi)
16      };
```

The iterative refinement seeks updates that are polynomial in the sub-patch field. Here gather the coefficients for later reference, where `maxn` is the maximum order of the polynomial: increase for higher-order or more nonlinear problems.

```
17 maxn:=4;
18 operator cl,cc,cr,cb;
19 cvars:=cg.(for n:=0:maxn join {cl(n),cc(n),cr(n),cb(n)})$
20 czero:=part(solve(cvars,cvars),1)$
```

Use various differences of macroscale grid values for the interpolative coupling. `uud(p)` is either $\delta^p U$ or $\mu\delta^p U$ when $p$ is even or odd respectively. We could analyse a low-order interpolation between patches to high-order in coupling parameter $\gamma$ in order to demonstrate convergence (hopefully).

```
21 maxint:=floor((maxn+1)/2);
22 array uud(2*maxint);
23 uud(0):=uu(j)$ % core value
24 uud(1):=(uu(j+1)-uu(j-1))/2$ % mu*delta
25 for p:=2:2*maxint do uud(p):= % delta^2 of two orders lower
26      sub(j=j+1,uud(p-2))-2*uud(p-2)+sub(j=j-1,uud(p-2))$
```

Recalling the shift operator $\mathcal{E}^{\pm r/2} = (1 \pm \mu\delta + \frac{1}{2}\delta^2)^{r/2}$ so expand this form in a Taylor series and then evaluate.

```
27 f:=taylor((1+eps)^xi,eps,0,maxint)$
28 f:=taylortostandard(f);
29 epxi:=(sub(eps=+mu*delta+delta^2/2,f)
30      where mu^2=>1+delta^2/4)$
31 emxi:=(sub(eps=-mu*delta+delta^2/2,f)
32      where mu^2=>1+delta^2/4)$
```

```
33 intr:=sub(xi=r/2,epxi);
34 intl:=sub(xi=r/2,emxi);
```
Then express these formulas as interpolation of the macroscale grid values $U_j$.
```
35 let gam^2=>gamma*epsilon;
36 uuintr:=uud(0)+sub(delta=uud(1)*gam,(-1+intr
37     where {delta^~~p=>uud(p)*gam^p,mu=>gam}));
38 uuintl:=uud(0)+sub(delta=uud(1)*gam,(-1+intl
39     where {delta^~~p=>uud(p)*gam^p,mu=>gam}));
```
Loop to seek corrections until residuals are smaller than specified order of error. This particular algorithm uses a single order parameter, $\epsilon = $ `epsilon`, for all small parameters (as in above multiplication of `gamma`): it assumes that the solution at each and every order in $\epsilon$ is found correctly and so does not need any subsequent refinement. The algorithm proceeds to higher and higher order in $\epsilon$, truncating to the specified orders of error in the parameters, until all residuals are zero. Because we do not truncate in $\epsilon$ then the residual calculation is exact to the specified order in parameters.
```
40 let { gamma^2=>0 , alpha^2=>0 };
41 for it:=1:99 do begin
```
Here we add a general form to the unknown fields. These coefficients are to be chosen to eliminate residuals in $\epsilon^{\tt it}$.
```
42 g:=g+cg*epsilon^it;
43 ul:=ul+epsilon^it*(for n:=0:maxn sum cl(n)*xi^n);
44 uc:=uc+epsilon^it*(for n:=0:maxn sum cc(n)*xi^n);
45 ur:=ur+epsilon^it*(for n:=0:maxn sum cr(n)*xi^n);
46 ub:=ub+epsilon^it*(for n:=0:maxn sum cb(n)*xi^n);
```
Compute the thirteen residuals of the advection-diffusion PDE in the four regions, the $C^1$ continuity conditions, and the patch amplitude condition. In the control use the action regions averages, and their difference with interpolation of core-averages from neighbouring patches. The control may be easiest to see in terms of its reciprocal, `rmu` $= 1/\mu$.
```
47 mu:=1/rmu;
48 ampl:=meanc(uc,xi)-uu(j);
49 pdel:=-df(ul,t)+kk*df(ul,x,x)-alpha*epsilon*df(ul,x)
50     +mu*kk/r^2/hh^2*(uuintl-meanl(ul,xi));
51 pdec:=-df(uc,t)+kk*df(uc,x,x)-alpha*epsilon*df(uc,x);
52 pder:=-df(ur,t)+kk*df(ur,x,x)-alpha*epsilon*df(ur,x)
53     +mu*kk/r^2/hh^2*(uuintr-meanr(ur,xi));
54 pdeb:=-df(ub,t)+kk*df(ub,x,x)-alpha*epsilon*df(ub,x);
55 c0bl:=sub(xi=-3*qr,ul)-sub(xi=5*qr,ub);
56 c0lc:=sub(xi=-qr,uc-ul);
57 c0cr:=sub(xi=+qr,ur-uc);
```

```
58 c0rb:=sub(xi=+3*qr,ub-ur);
59 c1bl:=sub(xi=-3*qr,df(ul,x))-sub(xi=5*qr,df(ub,x));
60 c1lc:=sub(xi=-qr,df(uc,x)-df(ul,x));
61 c1cr:=sub(xi=+qr,df(ur,x)-df(uc,x));
62 c1rb:=sub(xi=+3*qr,df(ub,x)-df(ur,x));
63 ress:={ampl, pdel, pdec, pder, pdeb, c0bl,
64     c0lc, c0cr, c0rb, c1bl, c1lc, c1cr, c1rb};
```
Monitor the progress of the residuals.
```
65 write lengthress:=map(length(~a),sub(czero,ress));
```
Solve for updates of order $\epsilon^{it}$ to the field and evolution.
```
66 eqns:=(foreach res in ress join
67     coeff(coeffn(res,epsilon,it),xi));
68 write lengtheqnsvars:={length(eqns),length(cvars)};
69 soln:=solve(eqns,cvars)$
70 if length(soln)=0 then rederr("******* no solution found");
71 g:=sub(soln,g);
72 ul:=sub(soln,ul);
73 uc:=sub(soln,uc);
74 ur:=sub(soln,ur);
75 ub:=sub(soln,ub);
```
Exit the loop when all residuals are zero.
```
76 showtime;
77 if sub(czero,ress)={0,0,0,0,0,0,0,0,0,0,0,0,0}
78    then write it:=it+100000;
79 end;
```
The single ordering parameter $\epsilon$ is no longer needed.
```
80 epsilon:=1$
```
Find the equivalent PDE to the discretisation where `du(n)` denotes the $n$th spatial derivative.
```
81 operator du;
82 factor hh,r;
83 erro:=deg((1+gamma)^9,gamma)*2+1;
84 equivpde:=(g where uu(~k)=>(du(0)+for m:=1:erro+2 sum
85                 ((k-j)*hh)^m*du(m)/factorial(m)) )$
86 for m:=erro:erro do let hh^m=>0;
87 equivpde:=equivpde;
```
Draw an example of the subpatch fields.
```
88 load_package gnuplot;
89 r:=1; rmu:=1/30; kk:=hh:=gamma:=1; alpha:=1;
90 procedure u(x);
91   (if x>3*qr then sub(xi=x,ub)
```

```
92  else if x>qr then sub(xi=x,ur)
93  else if x>-qr then sub(xi=x,uc)
94  else sub(xi=x,ul));
95 operator myu;
96 let myu(~x,~p)=> coeffn(u(x),uu(j+p),1) when numberp x;
97 plot(myu(x,0),myu(x,1),myu(x,-1),x=(-3*qr .. 5*qr)
98     ,size="0.7,0.5",title="",xlabel="",ylabel=""
99     ,terminal="postscript color",output="mmdcmppsX.eps");
```
Fin.
```
100 end;
```